\documentclass[10pt]{article}
\usepackage[plain]{fullpage}
\usepackage{times,amsmath,amsthm,amssymb,graphicx,xspace,epsfig,xcolor}
\usepackage{algorithm}\usepackage{algorithmic}
\usepackage{pgf}
\usepackage{pgffor}
\usepackage{xspace}
\usepackage{color}
\usepackage{comment}
\usepackage{authblk}
\usepackage{multirow}
\usepackage{hyperref}
\usepackage{cleveref}
\usepackage{dsfont}
\usepackage{caption}
\usepackage{subcaption}
\usepackage{soul}
\usepackage{tikz, tkz-graph, tkz-berge}
\usetikzlibrary{fit, backgrounds, matrix, arrows.meta}
\usetikzlibrary{calc,arrows,positioning}
\usetikzlibrary{decorations}
\usetikzlibrary{decorations.pathmorphing}
\usetikzlibrary{decorations.pathreplacing}
\usetikzlibrary{decorations.shapes}
\usetikzlibrary{decorations.text}
\usetikzlibrary{decorations.markings}
\usetikzlibrary{decorations.fractals}
\usetikzlibrary{decorations.footprints}
\usetikzlibrary{patterns}

\newtheorem{theorem}{Theorem}
\newtheorem{corollary}[theorem]{Corollary}

\newtheorem{lemma}[theorem]{Lemma}
\theoremstyle{definition}
\newtheorem{remark}[theorem]{Remark}
\newtheorem{question}[theorem]{Question}

\newtheorem{conjecture}[theorem]{Conjecture}

\newcommand{\ind}[1]{\langle #1 \rangle}

\newcommand{\defproblem}[3]{
     \vspace{1mm}
    \noindent\fbox{
     \begin{minipage}{0.96\textwidth}
     \begin{tabular*}{\textwidth}{@{\extracolsep{\fill}}lr} #1  \\ \end{tabular*}
     {\bf{Input:}} #2 \\
     {\bf{Question:}} #3
     \end{minipage}
     }
     \vspace{1mm}
}

\newcounter{claim}
\counterwithin{claim}{theorem}

\newenvironment{claim}[1][]
{\refstepcounter{claim}\vspace{1ex}\noindent{\bf Claim~\theclaim : }\it}{\vspace{1ex}}

\newenvironment{proofclaim}[1][]
	{\par\noindent {\it Proof of claim}. }{ \hfill$\lozenge$\par\vspace{11pt}}

\newenvironment{subproof}{\par\noindent {\it Proof of claim}.\ }{\hfill$\lozenge$\par\vspace{11pt}}

\DeclareMathOperator{\diff}{diff}

\definecolor{g-green}{rgb}{0.235, 0.659, 0.322}

\newcommand{\dic}{\vec{\chi}}
\newcommand{\bid}{\overleftrightarrow}

\begin{document}

\title{Strengthening the Directed Brooks' Theorem for oriented graphs and consequences on digraph redicolouring \thanks{Research supported by research grant
    DIGRAPHS ANR-19-CE48-0013 and by the French government, through the EUR DS4H Investments in the Future project managed by the National Research Agency (ANR) with the reference number ANR-17-EURE-0004.}}

\author{Lucas Picasarri-Arrieta}
\date{}

\maketitle
\vspace{-10mm}
\begin{center}
{\small 
Universit\'e C\^ote d'Azur, CNRS, I3S, INRIA, Sophia Antipolis, France\\
\texttt{lucas.picasarri-arrieta@inria.fr}\\
}
\end{center}

\maketitle

\begin{abstract}
    Let $D=(V,A)$ be a digraph. We define $\Delta_{\max}(D)$ as the maximum of $\{ \max(d^+(v),d^-(v)) \mid v \in V \}$ and  $\Delta_{\min}(D)$ as the maximum of $\{ \min(d^+(v),d^-(v)) \mid v \in V \}$. It is known that the dichromatic number of $D$ is at most $\Delta_{\min}(D) + 1$.
    In this work, we prove that every digraph $D$ which has dichromatic number exactly $\Delta_{\min}(D) + 1$ must contain the directed join of $\overleftrightarrow{K_r}$ and  $\overleftrightarrow{K_s}$ for some $r,s$ such that $r+s = \Delta_{\min}(D) + 1$, except if $\Delta_{\min}(D) = 2$ in which case $D$ must contain a digon. In particular, every oriented graph $\vec{G}$ with $\Delta_{\min}(\vec{G}) \geq 2$ has dichromatic number at most $\Delta_{\min}(\vec{G})$.
    
    Let $\vec{G}$ be an oriented graph of order $n$ such that $\Delta_{\min}(\vec{G}) \leq 1$. Given two 2-dicolourings of $\vec{G}$, we show that we can transform one into the other in at most $n$ steps, by recolouring one vertex at each step while maintaining a dicolouring at any step. Furthermore, we prove that, for every oriented graph $\vec{G}$ on $n$ vertices, the distance between two $k$-dicolourings is at most $2\Delta_{\min}(\vec{G})n$ when $k\geq \Delta_{\min}(\vec{G}) + 1$.
    
    We then extend a theorem of Feghali, Johnson and Paulusma to digraphs. We prove that, for every digraph $D$ with $\Delta_{\max}(D) = \Delta \geq 3$ and every $k\geq \Delta +1$, the $k$-dicolouring graph of $D$ consists of isolated vertices and at most one further component that has diameter at most $c_{\Delta}n^2$, where $c_{\Delta} = O(\Delta^2)$ is a constant depending only on $\Delta$.
\end{abstract}

\bibliographystyle{plain}

\section{Introduction}
\label{section:introduction}

\subsection{Graph (re)colouring}

 Given a graph $G=(V,E)$, a \textit{$k$-colouring} of $G$ is a function $c: V \xrightarrow{} \{1,\dots,k\}$ such that, for every edge $xy\in E$, we have $c(x) \neq c(y)$. So for every $i\in \{1,\dots,k\}$, $c^{-1}(i)$ induces an independent set on $G$. The \textit{chromatic number} of $G$, denoted by $\chi(G)$, is the smallest $k$ such that $G$ admits a $k$-colouring. The \textit{maximum degree} of $G$, denoted by $\Delta(G)$, is the degree of the vertex with the greatest number of edges incident to it.  A simple greedy procedure shows that, for any graph $G$, $\chi (G) \leq \Delta(G) + 1$. The celebrated theorem of Brooks~\cite{brooksMPCPS37} characterizes the graphs for which equality holds.
 
\begin{theorem}[Brooks, \cite{brooksMPCPS37}]\label{thm:brooks}
    A connected graph $G$ satisfies $\chi (G) = \Delta(G) + 1$ if and only if $G$ is an odd cycle or a complete graph.
\end{theorem}

 For any $k\geq \chi(G)$, the \textit{$k$-colouring graph} of $G$, denoted by ${\cal C}_k(G)$, is the graph whose vertices are the $k$-colourings of $G$ and in which two $k$-colourings are adjacent if they differ by the colour of exactly one vertex.
 A path between two given colourings in ${\cal C}_k(G)$ corresponds to a {\it recolouring sequence}, that is a sequence of pairs composed of a vertex of $G$, which is going to receive a new colour, and a new colour for this vertex.
 If ${\cal C}_k(G)$ is connected, we say that $G$ is \textit{$k$-mixing}. A $k$-colouring of $G$ is \textit{$k$-frozen} if it is an isolated vertex in ${\cal C}_k(G)$. The graph $G$ is \textit{$k$-freezable} if it admits a $k$-frozen colouring. In the last fifteen years, since the papers of Cereceda, van den Heuvel and Johnson~\cite{cerecedaJGT67,cerecedaEJC30}, graph recolouring has been studied by many researchers in graph theory.
 We refer the reader to the PhD thesis of Bartier~\cite{bartierTHESIS} for a complete overview on graph recolouring and to the surveys of van Heuvel~\cite{heuvel13} and Nishimura~\cite{Nishimura18} for reconfiguration problems in general.
 Feghali, Johnson and Paulusma~\cite{feghaliJGT83} proved the following analogue of Brooks' Theorem for graphs recolouring.

\begin{theorem}[\cite{feghaliJGT83}]\label{thm:feghali}
    Let $G=(V,E)$ be a connected graph with $\Delta(G) = \Delta \geq 3$, $k\geq \Delta+1$, and $\alpha$, $\beta$ two $k$-colourings of $G$. Then at least one of the following holds:
    \begin{itemize}
        \item $\alpha$ is $k$-frozen, or
        \item $\beta$ is $k$-frozen, or
        \item there is a recolouring sequence of length at most $c_{\Delta}|V|^2$ between $\alpha$ and $\beta$, where $c_\Delta = O(\Delta)$ is a constant depending on $\Delta$.
    \end{itemize}
\end{theorem}

Considering graphs of bounded maximum degree, Theorem~\ref{thm:feghali} has been very recently improved by Bousquet, Feuilloley, Heinrich and Rabie, who proved the following.
\begin{theorem}[\cite{bousquetARXIV220308885}]
    \label{thm:bousquet_al}
 Let $G=(V,E)$ be a connected graph with $\Delta(G) = \Delta \geq 3$, $k\geq \Delta+1$, and $\alpha$, $\beta$ two $k$-colourings of $G$. Then at least one of the following holds:
    \begin{itemize}
        \item $\alpha$ is $k$-frozen, or
        \item $\beta$ is $k$-frozen, or
        \item there is a recolouring sequence of length at most $O(\Delta^{c\Delta}|V|)$ between $\alpha$ and $\beta$, where $c$ is a constant.
    \end{itemize}
\end{theorem}
 
\subsection{Digraph (re)dicolouring}

In this paper, we are looking for extensions of the previous results on graphs colouring and recolouring to digraphs.

Let $D$ be a digraph. A \textit{digon} is a pair of arcs in opposite directions between the same vertices. A \textit{simple arc} is an arc which is not in a digon. For any two vertices $x,y\in V(D)$, the digon $\{xy,yx\}$ is denoted by $[x,y]$. The \textit{digon graph} of $D$ is the undirected graph with vertex set $V(D)$ in which $uv$ is an edge if and only if $[u,v]$ is a digon of $D$. 
An \textit{oriented graph} is a digraph with no digon. The \textit{bidirected graph} associated to a graph $G$, denoted by $\bid{G}$,  is the digraph obtained from $G$, by replacing every edge by a digon. The \textit{underlying graph} of $D$, denoted by $UG(D)$, is the undirected graph $G$ with vertex set $V(D)$ in which $uv$ is an edge if and only if $uv$ or $vu$ is an arc of $D$.

Let $v$ be a vertex of a digraph $D$. The {\it out-degree} (resp. {\it in-degree}) of $v$, denoted by $d^+(v)$ (resp. $d^-(v)$), is the number of arcs leaving (resp. entering) $v$. We define the \textit{maximum degree} of $v$ as $d_{\max}(v) = \max\{d^+(v), d^-(v)\}$, and the \textit{minimum degree} of $v$ as $d_{\min}(v) = \min\{d^+(v), d^-(v)\}$. We can then define the corresponding maximum degrees of $D$: $\Delta_{\max}(D) = \max_{v\in V(D)} (d_{\max}(v))$ and $\Delta_{\min}(D) = \max_{v\in V(D)} (d_{\min}(v))$. A digraph $D$ is {\it $\Delta$-diregular} if, for every vertex $v\in V(D)$, $d^-(v) = d^+(v) = \Delta$.

In 1982, Neumann-Lara~\cite{neumannlaraJCT33} introduced the notions of dicolouring and dichromatic number, which generalize the ones of colouring and chromatic number. 
A \textit{$k$-dicolouring} of $D$ is a function $c: V(D) \rightarrow{} \{1,\dots,k\}$ such that $c^{-1}(i)$ induces an acyclic subdigraph in $D$ for each $i \in \{1,\dots,k\}$. The \textit{dichromatic number} of $D$, denoted by $\dic(D)$, is the smallest $k$ such that $D$ admits a $k$-dicolouring. 
There is a one-to-one correspondence between the $k$-colourings of a graph $G$ and the $k$-dicolourings of the associated bidirected graph $\bid{G}$, and in particular $\chi(G) = \dic(\bid{G})$. Hence every result on graph colourings can be seen as a result on dicolourings of bidirected graphs, and it is natural to study whether the result can be extended to all digraphs.

The directed version of Brooks' Theorem was first proved by Harutyunyan and Mohar in~\cite{harutyunyanJDM25} (see also~\cite{aboulkerDM113193}). Aboulker and Aubian gave four new proofs of the following theorem in~\cite{aboulkerDM113193}.
\begin{theorem}[{\sc Directed Brooks' Theorem}]\label{thm:brooksdirected}
    Let $D$ be a connected digraph. Then $\dic(D) \leq \Delta_{\max}(D)+1$ and equality holds if and only if one of the following occurs:
    \begin{itemize}
        \item $D$ is a directed cycle, or
        \item $D$ is a bidirected odd cycle, or
        \item $D$ is a bidirected complete graph (of order at least $4$).
    \end{itemize}
\end{theorem}

It is easy to prove, by a simple greedy procedure, that every digraph $D$ can be dicoloured with $\Delta_{\min}(D)+1$ colours. Hence, one can wonder if Brooks' Theorem can be extended to digraphs using $\Delta_{\min}(D)$ instead of $\Delta_{\max}(D)$. Unfortunately, Aboulker and Aubian~\cite{aboulkerDM113193} proved that, given a digraph $D$, deciding whether $D$ is $\Delta_{\min}(D)$-dicolourable is NP-complete. Thus, unless P=NP, we cannot expect an easy characterization of digraphs satisfying $\dic(D) = \Delta_{\min}(D) + 1$. 

Let the {\it maximum geometric mean} of a digraph $D$ be $\Tilde{\Delta}(D) = \max \{ \sqrt{d^+(v)d^-(v)} \mid v \in V(D) \}$. By definition we have $\Delta_{\min}(D) \leq \Tilde{\Delta}(D) \leq \Delta_{\max}(D)$. Restricted to oriented graphs, Harutyunyan and Mohar~\cite{harutyunyanEJC18} have strengthened Theorem~\ref{thm:brooksdirected} by proving the following.
\begin{theorem}[Harutyunyan and Mohar~\cite{harutyunyanEJC18}]
    There is an absolute constant $\Delta_1$ such that every oriented graph $\vec{G}$ with $\Tilde{\Delta}(\vec{G}) \geq \Delta_1$ has $\dic(\vec{G}) \leq (1-e^{-13})\Tilde{\Delta}(\vec{G})$. 
\end{theorem}

In Section~\ref{section:brooks}, we give another strengthening of Theorem~\ref{thm:brooksdirected} on a large class of digraphs which contains oriented graphs. The \textit{directed join} of $H_1$ and $H_2$, denoted by $H_1 \Rightarrow H_2$, is the digraph obtained from disjoint copies of $H_1$ and $H_2$ by adding all arcs from the copy of $H_1$ to the copy of $H_2$ ($H_1$ or $H_2$ may be empty). 

\begin{theorem}\label{thm:brooks-deltamin}
Let $D$ be a digraph.
If $D$ is not $\Delta_{\min}(D)$-dicolourable, then one of the following holds: 
\begin{itemize}
    \item $\Delta_{\min}(D) \leq 1$, or
    \item $\Delta_{\min}(D) = 2$ and $D$ contains $\bid{K_2}$, or
    \item $\Delta_{\min}(D) \geq 3$ and $D$ contains $\bid{K_r} \Rightarrow \bid{K_s}$, for some $r,s \geq 0$ such that $r+s = \Delta_{\min}(D) + 1$.
\end{itemize}
\end{theorem}

In particular, the following is a direct consequence of Theorem~\ref{thm:brooks-deltamin}.

\begin{corollary}
    \label{cor:brooks-deltamin}
    Let $D$ be a digraph. If $\dic(D) = \Delta_{\min}(D) + 1$, then $D$ contains the complete bidirected graph on $\left \lceil \frac{\Delta_{\min}+1}{2} \right \rceil$ vertices as a subdigraph.
\end{corollary}

Corollary~\ref{cor:brooks-deltamin} is best possible: if we restrict $D$ to not contain the complete bidirected graph on $\left \lceil \frac{\Delta_{\min}+1}{2} \right \rceil + 1$, then we show that deciding whether $\dic(D) \leq \Delta_{\min}(D)$ remains NP-complete (Theorem~\ref{thm:np_complete_deltamin}). 
Moreover, since an oriented graph does not contain any digon, Corollary~\ref{cor:brooks-deltamin} implies the following:

\begin{corollary}
    \label{cor:brooks-oriented}
    Let $\vec{G}$ be an oriented graph. If $\Delta_{\min}(\vec{G}) \geq 2$, then $\dic(\vec{G}) \leq \Delta_{\min}(\vec{G})$.
\end{corollary}

\medskip

For any $k\geq \dic(D)$, the \textit{$k$-dicolouring graph} of $D$, denoted by ${\cal D}_k(D)$, is the graph whose vertices are the $k$-dicolourings of $D$ and in which two $k$-dicolourings are adjacent if they differ by the colour of exactly one vertex.  Observe that ${\cal C}_k(G) = {\cal D}_k(\bid{G})$ for any bidirected graph $\bid{G}$.  A \textit{redicolouring sequence} between two dicolourings is a path between these dicolourings in ${\cal D}_k(D)$.  The digraph $D$ is \textit{$k$-mixing} if ${\cal D}_k(D)$ is connected. A $k$-dicolouring of $D$ is \textit{$k$-frozen} if it is an isolated vertex in ${\cal D}_k(D)$. The digraph $D$ is \textit{$k$-freezable} if it admits a $k$-frozen dicolouring. A vertex $v$ is \textit{blocked} to its colour in a dicolouring $\alpha$ if, for every colour $c\neq \alpha(v)$, recolouring $v$ to $c$ in $\alpha$ creates a monochromatic directed cycle. 

Digraph redicolouring was first introduced in~\cite{papierAMADEUS}, where the authors generalized different results on graph recolouring to digraphs, and proved some specific results on oriented graphs redicolouring. In particular, they studied the $k$-dicolouring graph of digraphs with bounded degeneracy or bounded maximum average degree, and they show that finding a redicolouring sequence between two given $k$-dicolourings of a digraph is PSPACE-complete.
Dealing with the maximum degree of a digraph, they proved that, given an orientation of a subcubic graph $\vec{G}$ on $n$ vertices, its $2$-dicolouring graph ${\cal D}_2(\vec{G})$ is connected and has diameter at most $2n$ and they asked if this bound can be improved. We answer this question in Section~\ref{section:recolouring_oriented} by proving the following theorem.
\begin{theorem}\label{thm:mindegree-1}
Let $\vec{G}$ be an oriented graph of order $n$ such that $\Delta_{\min}(\vec{G}) \leq 1$. Then ${\cal D}_2(\vec{G})$ is connected and has diameter exactly $n$.
\end{theorem}

In particular, if $\vec{G}$ is an orientation  of a subcubic graph, then $\Delta_{\min}(\vec{G}) \leq 1$ (because $d^+(v) + d^-(v) \leq 3$ for every vertex $v$), and so ${\cal D}_2(\vec{G})$ has diameter exactly $n$. Furthermore, we prove the following as a consequence of Corollary~\ref{cor:brooks-oriented} and Theorem~\ref{thm:mindegree-1}.

\begin{corollary}\label{cor:Deltamin}
    Let $\vec{G}$ be oriented graph of order $n$ with $\Delta_{\min}(\vec{G}) = \Delta \geq 1$, and let $k \geq \Delta + 1$. Then ${\cal D}_k(\vec{G})$ is connected and has diameter at most $2\Delta n$.
\end{corollary}

Corollary~\ref{cor:Deltamin} does not hold for digraphs in general: indeed, $\bid{P_n}$, the bidirected path on $n$ vertices, satisfies $\Delta_{\min}(\bid{P_n}) = 2$ and ${\cal D}_3(\bid{P_n}) = {\cal C}_3(P_n)$ has diameter $\Omega(n^2)$, as proved in~\cite{bonamyJCO27}. 

Finally in Section~\ref{section:recolouring_directed}, we extend Theorem~\ref{thm:feghali} to digraphs. 

\begin{theorem}\label{thm:extension_feghali}
    Let $D=(V,A)$ be a connected digraph with $\Delta_{\max}(D) = \Delta \geq 3$, $k\geq \Delta+1$, and $\alpha$, $\beta$ two $k$-dicolourings of $D$. Then at least one of the following holds:
    \begin{itemize}
        \item $\alpha$ is $k$-frozen, or
        \item $\beta$ is $k$-frozen, or
        \item there is a redicolouring sequence of length at most $c_{\Delta}|V|^2$ between $\alpha$ and $\beta$,
    where $c_{\Delta} = O(\Delta^2)$ is a constant depending only on $\Delta$.
    \end{itemize}
\end{theorem}

Furthermore, we prove that a digraph $D$ is $k$-freezable only if $D$ is bidirected and its underlying graph is $k$-freezable. Thus, an obstruction in Theorem~\ref{thm:extension_feghali} is exactly the bidirected graph of an obstruction in Theorem~\ref{thm:feghali}.

\section{Strengthening of Directed Brooks' Theorem for oriented graphs}
\label{section:brooks}
A digraph $D$ is \textit{$k$-dicritical} if $\dic(D) = k$ and for every vertex $v\in V(D)$, $\dic (D-v) < k$. Observe that every digraph with dichromatic number at least $k$ contains a $k$-dicritical subdigraph. 

Let $\mathcal{F}_2$ be $\{ \overleftrightarrow{K_2} \}$, and for each $\Delta\geq 3$, we define $\mathcal{F}_\Delta = \{ \bid{K_r} \Rightarrow \bid{K_s} \mid r,s\geq 0~\mbox{and}~r+s = \Delta+1 \}$. A digraph $D$ is \textit{$\mathcal{F}_\Delta$-free} if it does not contain $F$ as a subdigraph, for any $F\in \mathcal{F}_\Delta$. Theorem~\ref{thm:brooks-deltamin} can then be reformulated as follows:

\begingroup
    \def\thetheorem{\ref{thm:brooks-deltamin}}
        \begin{theorem}
            Let $D$ be a digraph with $\Delta_{\min}(D) = \Delta \geq 2$. If $D$ is $\mathcal{F}_\Delta$-free, then $\dic(D) \leq \Delta$.
        \end{theorem}
    \addtocounter{theorem}{-6}
\endgroup

\begin{proof}
    Let $D$ be a digraph such that $ \Delta_{\min}(D) = \Delta \geq 2$ and $\dic(D) = \Delta + 1$. We will show that $D$ contains some $F\in \mathcal{F}_\Delta$ as a subdigraph. 
    
    Let $(X,Y)$ be a partition of $V(D)$ such that for each $x\in X$, $d^+(x)\leq \Delta$, and for each $y\in Y$, $d^-(y) \leq \Delta$. We define the digraph $\Tilde{D}$ as follows:
    \begin{itemize}
        \item $V(\Tilde{D}) = V(D)$,
        \item $A(\Tilde{D}) = A(D\ind{X}) \cup A(D\ind{Y}) \cup \{xy,yx \mid xy\in A(D), x\in X, y\in Y\}$.
    \end{itemize}
    
    \begin{claim}\label{claim:imply-2dic}
        $\dic(\Tilde{D}) \geq \Delta+1$.
    \end{claim}
    \begin{subproof}
        Assume for a contradiction that there exists a $\Delta$-dicolouring $c$ of $\Tilde{D}$. Then $D$, coloured with $c$, must contain a monochromatic directed cycle $C$. Now $C$ is not contained in $X$ nor $Y$, for otherwise $C$ would be a monochromatic directed cycle of $D\ind{X}$ or $D\ind{Y}$ and so a monochromatic directed cycle of $\Tilde{D}$. Thus $C$ contains an arc $xy$ from $X$ to $Y$. But then, $[x,y]$ is a monochromatic digon in $\Tilde{D}$, a contradiction.
    \end{subproof}
    
    Since $\dic(\Tilde{D}) \geq \Delta+1$, there is a ($\Delta+1$)-dicritical subdigraph $H$ of $\Tilde{D}$. By dicriticality of $H$, for every vertex $v\in V(H)$, $d_H^+(v)\geq \Delta$ and $d_H^-(v) \geq \Delta$, for otherwise a $\Delta$-dicolouring of $H-v$ could be extended to $H$ by choosing for $v$ a colour which is not appearing in its out-neighbourhood or in its in-neighbourhood. We define $X_H$ as $X\cap V(H)$ and $Y_H$ as $Y\cap V(H)$.  Note that both $H\ind{X_H}$ and $H\ind{Y_H}$ are subdigraphs of $D$.
    
    \begin{claim}
        $H$ is $\Delta$-diregular.
    \end{claim}
    \begin{proofclaim}
        Let $\ell$ be the number of digons between $X_H$ and $Y_H$ in $H$.
        Observe that, by definition of $X$ and $H$, for each vertex $x\in X_H$, $d_H^+(x) = \Delta$. Note also that, in $H$, $\ell$ is exactly the number of arcs leaving $X_H$ and exactly the number of arcs entering $X_H$. We get:
        \begin{align*}
            \Delta|X_H| &= \sum_{x\in X_H}d^+_{H}(x)\\
                &=  \ell + |A(H\ind{X_H})|\\
                &=  \sum_{x\in X_H}d^-_{H}(x)
        \end{align*}
        which implies, since $H$ is dicritical, $d_H^+(x) = d_H^-(x) = \Delta$ for every vertex $x\in X_H$. Using a symmetric argument, we prove that $\Delta|Y_H| = \sum_{y\in Y_H}d^+_{H}(y)$, implying $d_H^+(y) = d_H^-(y) = \Delta$ for every vertex $y\in Y_H$.
    \end{proofclaim}
    
    Since $H$ is $\Delta$-diregular, then in particular $\Delta_{\max}(H) = \Delta$. Hence, because $\dic(H) = \Delta+1$, by Theorem~\ref{thm:brooksdirected}, either $\Delta=2$ and $H$ is a bidirected odd cycle, or $\Delta \geq 3$ and $H$ is the bidirected complete graph on $\Delta+1$ vertices.
    \begin{itemize}
        \item If $\Delta=2$ and $H$ is a bidirected odd cycle, then at least one digon of $H$ belongs to $H\ind{X_H}$ or $H\ind{Y_H}$, for otherwise $H$ would be bipartite (with bipartition $(X_H,Y_H)$). Since both $H\ind{X_H}$ and $H\ind{Y_H}$ are subdigraphs of $D$, this shows, as desired, that $D$ contains a copy of $\overleftrightarrow{K_2}$.
        \item If $k\geq 3$ and $H$ is the bidirected complete graph on $\Delta+1$ vertices, let $A_H$ be all the arcs from $Y_H$ to $X_H$. Then $D\ind{V(H)} \setminus A_H$ is a subdigraph of $D$ which belongs to $\mathcal{F}_\Delta$.
    \end{itemize}
 \vspace*{-12pt}   \end{proof}
\addtocounter{theorem}{+5}

Now we will justify that Corollary~\ref{cor:brooks-deltamin} is best possible. To do so, we prove that given a digraph $D$ which does not contain the bidirected complete graph on $\left \lceil \frac{\Delta_{\min}(D)+1}{2} \right \rceil + 1$ vertices, deciding if it is $\Delta_{\min}(D)$-dicolourable is NP-complete. 
We shall use a reduction from {\sc $k$-Dicolourability} which is defined as follows:\\
\defproblem{\textsc{$k$-Dicolourability}}{A digraph $D$}{Is $D$ $k$-dicolourable ?}
\smallskip

{\sc $k$-Dicolourability} is NP-complete for every fixed $k\geq 2$~\cite{bokalJGY46}. It remains NP-complete when we restrict to digraphs $D$ with $\Delta_{\min}(D) = k$~\cite{aboulkerDM113193}. 

\begin{theorem}
    \label{thm:np_complete_deltamin}
    For all $k\geq 2$, {\sc $k$-Dicolourability} remains NP-complete when restricted to digraphs $D$ satisfying $\Delta_{\min}(D) = k$ and not containing the bidirected complete graph on $\left \lceil \frac{k+1}{2} \right \rceil + 1$ vertices.
\end{theorem}
\begin{proof}
    Let $D=(V,A)$ be an instance of  {\sc $k$-Dicolourability} for some fixed $k\geq 2$. Then we build $D'=(V',A')$ as follows:
    \begin{itemize}
        \item For each vertex $x\in V$, we associate a copy of $S_x^- \Rightarrow  S_x^+$ where $S_x^-$ is the bidirected complete graph on $ \left \lfloor \frac{k+1}{2} \right \rfloor$ vertices, and $S_x^+$ is the bidirected complete graph on $ \left \lceil \frac{k+1}{2} \right \rceil$ vertices. 
        
        \item For each arc $xy\in A$, we associate all possible arcs $x^+y^-$ in $A'$, such that $x^+ \in S_x^+$ and $y^-\in S_y^-$.
    \end{itemize}
    First observe that $\Delta_{\min}(D') = k$. Let $v$ be a vertex of $D'$, if $v$ belongs to some $S_x^+$, then $d^-(v) = k$, otherwise it belongs to some $S_x^-$ and then $d^+(v) = k$.
    Then observe that $D'$ does not contain the bidirected complete graph on $\left \lceil \frac{k+1}{2} \right \rceil + 1$ vertices since every digon in $D'$ is contained in some $S_x^+$ or $S_x^-$.
    Thus we only have to prove that $\dic(D) \leq k$ if and only if $\dic(D') \leq k$ to get the result.
    
    \begin{itemize}
    \item Let us first prove that  $\dic(D) \leq k$ implies $\dic(D') \leq k$.\\
    Assume that $\dic(D) \leq k$. Let $\phi : V \xrightarrow{} \{1,\dots,k\}$ be a $k$-dicolouring of $D$. Let $\phi'$ be the $k$-dicolouring of $D'$ defined as follows: for each vertex $x\in V$, choose arbitrarily $x^- \in S_x^-$, $x^+\in S_x^+$, and set $\phi'(x^-) =\phi'(x^+) = \phi(x)$. Then choose a distinct colour for every other vertex $v$ in $S_x^- \cup S_x^+$, and set $\phi'(v)$ to this colour.
    We get that $\phi'$ must be a $k$-dicolouring of $D'$: for each $x\in V$, every vertex but $x^-$ in $S_x^-$ must be a sink in its colour class, and every vertex but $x^+$ in $S_x^+$ must be a source in its colour class. Thus if $D'$, coloured with $\phi'$, contains a monochromatic directed cycle $C'$, then $C'$ must be of the form $x_1^-x_1^+x_2^-x_2^+\cdots x_\ell^-x_\ell^+x_1^-$. But then $C=x_1x_2\cdots x_\ell x_1$ is a monochromatic directed cycle in $D$ coloured with $\phi$: a contradiction.
    
    \item  Reciprocally, let us prove that $\dic(D') \leq k$ implies $\dic(D) \leq k$.\\ Assume that $\dic(D') \leq k$.
    Let $\phi' : V' \xrightarrow{} \{1,\dots,k\}$ be a $k$-dicolouring of $D'$. Let $\phi$ be the $k$-dicolouring of $D$ defined as follows. For each vertex $x\in V$, we know that $|S_x^+\cup S_x^-| = k+1$, thus there must be two vertices $x^+$ and $x^-$ in $S_x^+\cup S_x^-$ such that $\phi'(x^+) = \phi'(x^-)$. Moreover, since both $S_x^+$ and $S_x^-$ are bidirected, one of these two vertices belongs to $S_x^+$ and the other one belongs to $S_x^-$. We assume without loss of generality $x^+ \in S_x^+$ and $x^- \in S_x^-$. Then we set $\phi(x) = \phi'(x^+)$. We get that $\phi$ must be a $k$-dicolouring of $D$. If $D$, coloured with $\phi$, contains a monochromatic directed cycle $C = x_1x_2\cdots x_\ell x_1$, then $C'=x_1^-x_1^+x_2^-x_2^+\cdots x_\ell^-x_\ell^+ x_1^-$ is a monochromatic directed cycle in $D'$ coloured with $\phi'$, a contradiction.
    \end{itemize}
 \vspace*{-12pt}   \end{proof}

\section{Redicolouring oriented graphs}
\label{section:recolouring_oriented}

In this section, we restrict to oriented graphs. We first prove Theorem~\ref{thm:mindegree-1}, let us restate it.

\begingroup
    \def\thetheorem{\ref{thm:mindegree-1}}
        \begin{theorem}
            Let $\vec{G}$ be an oriented graph of order $n$ such that $\Delta_{\min}(\vec{G}) \leq 1$. Then ${\cal D}_2(\vec{G})$ is connected and has diameter exactly $n$.
        \end{theorem}
    \addtocounter{theorem}{-1}
\endgroup

Observe that, if ${\cal D}_2(\vec{G})$ is connected, then its diameter must be at least $n$: for any 2-dicolouring $\alpha$, we can define its mirror $\bar{\alpha}$ where, for every vertex $v\in V(\vec{G})$, $\alpha(v) \neq \bar{\alpha}(v)$; then every redicolouring sequence between $\alpha$ and $\bar{\alpha}$ has length at least $n$.

\begin{lemma}\label{lemma:cycle}
    Let $C$ be a directed cycle of length at least $3$. Then ${\cal D}_2(C)$ is connected and has diameter exactly $n$.
\end{lemma}
\begin{proof}
    Let $\alpha$ and $\beta$ be any two 2-dicolourings of $C$. Let $x=\diff (\alpha,\beta) =  |\{v \in V(C) \mid \alpha (v) \neq \beta(v) \}|$. By induction on $x\geq 0$, let us show that there exists a path of length at most $x$ from $\alpha$ to $\beta$ in ${\cal D}_2(C)$. This clearly holds for $x = 0$ (i.e., $\alpha = \beta$). Assume $x > 0$ and the result holds for $x-1$. Let $v\in V(C)$ be  such that $\alpha(v) \neq \beta(v)$. 

    If $v$ can be recoloured in $\beta(v)$, then we recolour it and reach a new $2$-dicolouring $\alpha'$ such that $\diff(\alpha', \beta) = x - 1$ and the result holds by induction. Else if $v$ cannot be recoloured, then recolouring $v$ must create a monochromatic directed cycle, which must be $C$. Then there must be a vertex $v'$, different from $v$, such that $\beta(v) = \alpha(v') \neq \beta(v')$, and $v'$ can be recoloured. We recolour it and reach a new $2$-dicolouring $\alpha'$ such that $\diff(\alpha', \beta) = x - 1$ and the result holds by induction.
\end{proof}
\addtocounter{theorem}{-4}
We are now ready to prove Theorem~\ref{thm:mindegree-1}.
\begin{proof}[Proof of Theorem~\ref{thm:mindegree-1}]
Let $\alpha$ and $\beta$ be any two $2$-dicolourings of $\vec{G}$. We will show that there exists a redicolouring sequence of length at most $n$ between $\alpha$ and $\beta$. We may assume that $\vec{G}$ is strongly connected, otherwise we consider each strongly connected component independently. This implies in particular that $\vec{G}$ does not contain any sink nor source. Let $(X,Y)$ be a partition of $V(\vec{G})$ such that, for every $x\in X$, $d^+(x) = 1$, and for every $y\in Y$, $d^-(y) = 1$. 

Assume first that $\vec{G}\ind{X}$ contains a directed cycle $C$. Since every vertex in $X$ has exactly one out-neighbour, there is no arc leaving $C$. Thus, since $\vec{G}$ is strongly connected, $\vec{G}$ must be exactly $C$, and the result holds by Lemma~\ref{lemma:cycle}. Using a symmetric argument, we get the result when  $\vec{G}\ind{Y}$ contains a directed cycle.

Assume now that both $\vec{G}\ind{X}$ and $\vec{G}\ind{Y}$ are acyclic.  Thus, since every vertex in $X$ has exactly one out-neighbour, $\vec{G}\ind{X}$ is the union of disjoint and independent in-trees, that are oriented trees in which all arcs are directed towards the root.
We denote by $X_r$ the set of roots of these in-trees. Symmetrically, $\vec{G}\ind{Y}$ is the union of disjoint and independent out-trees (oriented trees in which all arcs are directed away from the root), and we denote by $Y_r$ the set of roots of these out-trees. Set $X_\ell=X \setminus X_r$ and $Y_\ell=Y \setminus Y_r$.
Observe that the arcs from $X$ to $Y$ form a perfect matching directed from $X_r$ to $Y_r$. We denote by $M_r$ this perfect matching. Observe also that there can be any arc from $Y$ to $X$.
Now we define $X_r^1$ and $Y_r^1$ two subsets of $X_r$ and $Y_r$ respectively, depending on the two $2$-dicolourings $\alpha$ and $\beta$, as follows:
\begin{align*}
    X_r^1 &= \{ x \mid xy \in M_r, \alpha(x) = \beta(y) \neq \alpha(y) = \beta(x) \}\\
    Y_r^1 &= \{ y \mid xy \in M_r, \alpha(x) = \beta(y) \neq \alpha(y) = \beta(x) \}
\end{align*}
Set $X_r^2 = X_r \setminus X_r^1$  and  $Y_r^2 = Y_r \setminus Y_r^1$. We denote by $M_r^1$ (respectively $M_r^2$) the perfect matching from $X_r^1$ to $Y_r^1$ (respectively from $X_r^2$ to $Y_r^2$). Figure~\ref{fig:partitioning-Delta-min-1} shows a partitioning of $V(\vec{G})$ into $X_r^1,X_r^2,X_\ell,Y_r^1,Y_r^2,Y_\ell$.

\begin{figure}[hbtp]
    \begin{minipage}{\linewidth}
        \begin{center}	
          \begin{tikzpicture}[thick,scale=1, every node/.style={transform shape}]
    	    \tikzset{vertex/.style = {circle,fill=black,minimum size=6pt, inner sep=0pt}}
    	    \tikzset{littlevertex/.style = {circle,fill=black,minimum size=0pt, inner sep=0pt}}
            \tikzset{edge/.style = {->,> = latex'}}
    	    
            \node[vertex,red] (1a) at  (0,0) {};
            \node[vertex] (2a) at  (0,-1) {};
            \node[vertex,red] (3a) at  (0,-2) {};
            \node[vertex,red] (4a) at  (0,-3) {};
            \node[vertex] (5a) at  (0,-4) {};
            \node[] (Xr1) at (-0.5,-0.5) {$X_r^1$};
            \node[] (Xr2) at (-0.5,-2.5) {$X_r^2$};
            \node[] (Xl) at (-2,-2) {$X_\ell$};
            \draw[dashed] (-2.5,0.5) rectangle ++(2.8,-5);
            \draw[dashed] (-1, 0.5) to (-1, -4.5);
            \draw[dashed] (-1, -1.5) to (0.3, -1.5);
            \foreach \i in {1,...,8}{
               \node [littlevertex]  (xl\i) at (-1.5,-0.58*\i+0.58) {};
               \node [littlevertex]  (xl2\i) at (-2,-0.58*\i+0.58) {};
            }
           \draw[edge]  (xl1) to (1a);
            \foreach \i in {2,...,3}{
               \draw[edge]  (xl\i) to (2a);
            }
            \foreach \i in {4,...,5}{
               \draw[edge]  (xl\i) to (3a);
            }
            \foreach \i in {6,...,8}{
               \draw[edge]  (xl\i) to (5a);
            }
            \foreach \i in {1,...,2}{
               \draw[edge]  (xl2\i) to (xl1);
            }
            \foreach \i in {5,...,6}{
               \draw[edge]  (xl2\i) to (xl5);
            }
            \foreach \i in {7,...,8}{
               \draw[edge]  (xl2\i) to (xl8);
            }
            \node[vertex] (1b) at  (1,0) {};
            \node[vertex,red] (2b) at  (1,-1) {};
            \node[vertex] (3b) at  (1,-2) {};
            \node[vertex,red] (4b) at  (1,-3) {};
            \node[vertex] (5b) at  (1,-4) {};
            \node[] (Yr1) at (1.5,-0.5) {$Y_r^1$};
            \node[] (Yr2) at (1.5,-2.5) {$Y_r^2$};
            \node[] (Yl) at (3,-2) {$Y_\ell$};
            \draw[dashed] (0.6,0.5) rectangle ++(2.9,-5);
            \draw[dashed] (2, 0.5) to (2, -4.5);
            \draw[dashed] (2, -1.5) to (0.7, -1.5);
            \foreach \i in {1,...,8}{
               \node [littlevertex]  (yl\i) at (2.5,-0.58*\i+0.58) {};
               \node [littlevertex]  (yl2\i) at (3,-0.58*\i+0.58) {};
            }
           \draw[edge]  (2b) to (yl3);
            \foreach \i in {1,...,2}{
               \draw[edge] (1b) to (yl\i);
            }
            \foreach \i in {4,...,4}{
               \draw[edge]  (3b) to (yl\i);
            }
            \foreach \i in {5,...,7}{
               \draw[edge] (4b) to (yl\i);
            }
            \foreach \i in {8,...,8}{
               \draw[edge] (5b) to (yl\i);
            }
            
            \foreach \i in {1,...,3}{
               \draw[edge]  (yl2) to (yl2\i);
            }
            \foreach \i in {5,...,6}{
               \draw[edge]  (yl5) to (yl2\i);
            }
            \foreach \i in {7,...,8}{
               \draw[edge] (yl8) to (yl2\i);
            }
            \foreach \i in {1,...,5}{
               \draw[edge] (\i a) to (\i b);
            }
            \draw [decorate,decoration={brace,amplitude=5pt,mirror,raise=4ex}] (-2.5,-4.3) -- (3.5,-4.3) node[midway,yshift=-3em]{$\vec{G}$ dicoloured with $\alpha$};
            
            \node[vertex] (1a_) at  (6.5,0) {};
            \node[vertex,red] (2a_) at  (6.5,-1) {};
            \node[vertex,red] (3a_) at  (6.5,-2) {};
            \node[vertex] (4a_) at  (6.5,-3) {};
            \node[vertex,red] (5a_) at  (6.5,-4) {};
            \node[] (Xr1_) at (6,-0.5) {$X_r^1$};
            \node[] (Xr2_) at (6,-2.5) {$X_r^2$};
            \node[] (Xl_) at (4.5,-2) {$X_\ell$};
            \draw[dashed] (4,0.5) rectangle ++(2.8,-5);
            \draw[dashed] (5.5, 0.5) to (5.5, -4.5);
            \draw[dashed] (5.5, -1.5) to (6.8, -1.5);
            \foreach \i in {1,...,8}{
               \node [littlevertex]  (xl_\i) at (5,-0.58*\i+0.58) {};
               \node [littlevertex]  (xl2_\i) at (4.5,-0.58*\i+0.58) {};
            }
           \draw[edge]  (xl_1) to (1a_);
            \foreach \i in {2,...,3}{
               \draw[edge]  (xl_\i) to (2a_);
            }
            \foreach \i in {4,...,5}{
               \draw[edge]  (xl_\i) to (3a_);
            }
            \foreach \i in {6,...,8}{
               \draw[edge]  (xl_\i) to (5a_);
            }
            \foreach \i in {1,...,2}{
               \draw[edge]  (xl2_\i) to (xl_1);
            }
            \foreach \i in {5,...,6}{
               \draw[edge]  (xl2_\i) to (xl_5);
            }
            \foreach \i in {7,...,8}{
               \draw[edge]  (xl2_\i) to (xl_8);
            }
            \node[vertex,red] (1b_) at  (7.5,0) {};
            \node[vertex] (2b_) at  (7.5,-1) {};
            \node[vertex] (3b_) at  (7.5,-2) {};
            \node[vertex,red] (4b_) at  (7.5,-3) {};
            \node[vertex,red] (5b_) at  (7.5,-4) {};
            \node[] (Yr1_) at (8,-0.5) {$Y_r^1$};
            \node[] (Yr2_) at (8,-2.5) {$Y_r^2$};
            \node[] (Yl_) at (9.5,-2) {$Y_\ell$};
            \draw[dashed] (7.1,0.5) rectangle ++(2.9,-5);
            \draw[dashed] (8.5, 0.5) to (8.5, -4.5);
            \draw[dashed] (8.5, -1.5) to (7.2, -1.5);
            \foreach \i in {1,...,8}{
               \node [littlevertex]  (yl_\i) at (9,-0.58*\i+0.58) {};
               \node [littlevertex]  (yl2_\i) at (9.5,-0.58*\i+0.58) {};
            }
           \draw[edge]  (2b_) to (yl_3);
            \foreach \i in {1,...,2}{
               \draw[edge] (1b_) to (yl_\i);
            }
            \foreach \i in {4,...,4}{
               \draw[edge]  (3b_) to (yl_\i);
            }
            \foreach \i in {5,...,7}{
               \draw[edge] (4b_) to (yl_\i);
            }
            \foreach \i in {8,...,8}{
               \draw[edge] (5b_) to (yl_\i);
            }
            \foreach \i in {1,...,3}{
               \draw[edge]  (yl_2) to (yl2_\i);
            }
            \foreach \i in {5,...,6}{
               \draw[edge]  (yl_5) to (yl2_\i);
            }
            \foreach \i in {7,...,8}{
               \draw[edge] (yl_8) to (yl2_\i);
            }
            \foreach \i in {1,...,5}{
               \draw[edge] (\i a_) to (\i b_);
            }
            \draw [decorate,decoration={brace,amplitude=5pt,mirror,raise=4ex}] (4,-4.3) -- (10,-4.3) node[midway,yshift=-3em]{$\vec{G}$ dicoloured with $\beta$};
          \end{tikzpicture}
      \caption{The partitioning of $V(\vec{G})$ into $X_r^1,X_r^2,X_\ell,Y_r^1,Y_r^2,Y_\ell$.}
      \label{fig:partitioning-Delta-min-1}
    \end{center}    
  \end{minipage}
\end{figure}
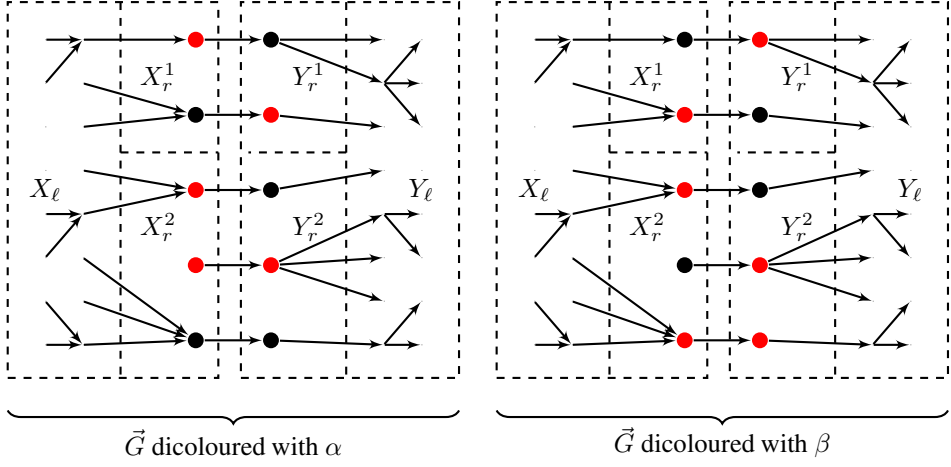
\begin{claim}
    There exists a redicolouring sequence of length $s_\alpha$ from $\alpha$ to some $2$-dicolouring $\alpha'$  and a redicolouring sequence of length $s_\beta$ from $\beta$ to some $2$-dicolouring $\beta'$ such that each of the following holds:
    \begin{itemize}
        \item[(i)] For any arc $xy \in M_r$, $\alpha'(x) \neq \alpha'(y)$ and $\beta'(x) \neq \beta'(y)$,
        \item[(ii)] For any arc $xy \in M^2_r$, $\alpha'(x) = \beta'(x)$ (and so $\alpha'(y) = \beta'(y)$ by (i)), and
        \item[(iii)] $s_\alpha + s_\beta \leq |X^2_r| + |Y^2_r|$.
    \end{itemize}
\end{claim}
\begin{proofclaim}
     We consider the arcs $xy$ of $M_r^2$ one after another and do the following recolourings  depending on the colours of $x$ and $y$ in both $\alpha$ and $\beta$ to get $\alpha'$ and $\beta'$.
    \begin{itemize}
        \item If $\alpha(x) = \alpha(y) = \beta(x) = \beta(y)$, then we recolour $x$ in both $\alpha$ and $\beta$;
        \item Else if $\alpha(x) = \alpha(y) \neq \beta(x) = \beta(y)$, then we recolour $x$ in $\alpha$ and we recolour $y$ in $\beta$;
        \item Else if $\alpha(x) = \beta(x) \neq \alpha(y) = \beta(y)$, then we do nothing;
        \item Else if $\alpha(x) \neq \alpha(y) = \beta(x) = \beta(y)$, then we recolour $x$ in $\beta$;
        \item Finally if $\alpha(y) \neq \alpha(x) = \beta(x) = \beta(y)$, then we recolour $y$ in $\beta$.
    \end{itemize}
    Each of these recolourings is valid because, when a vertex in $X_r^2$ (respectively $Y_r^2$) is recoloured, it gets a colour different from its only out-neighbour (respectively in-neighbour).
    Let $\alpha'$ and $\beta'$ be the the two resulting $2$-dicolourings. By construction, $\alpha'$ and $\beta'$ agree on $X_r^2\cup Y_r^2$. For each arc $xy \in M_r^2$, either $\alpha(x) = \alpha'(x)$ or $\alpha(y) = \alpha'(y)$, and the same holds for $\beta$ and $\beta'$. This implies that $s_\alpha + s_\beta \leq 2|M_r^2| = |X^2_r| + |Y^2_r|$.
\end{proofclaim}

\begin{claim}
There exists a redicolouring sequence from $\alpha'$ to some 2-dicolouring $\Tilde{\alpha}$ of length $s_\alpha'$ and a redicolouring sequence from $\beta'$ to some $2$-dicolouring $\Tilde{\beta}$ of length $s_\beta'$ such that each of the following holds:
    \begin{itemize}
        \item[(i)] $\Tilde{\alpha}$ and $\Tilde{\beta}$ agree on $V(\vec{G}) \setminus (X_r^1 \cup Y_r^1)$,
        \item[(ii)] $\alpha'$ and $\Tilde{\alpha}$ agree on $X_r \cup Y_r$,
        \item[(iii)] $\beta'$ and $\Tilde{\beta}$ agree on $X_r \cup Y_r$,
        \item[(iv)] $X_\ell \cup Y_\ell$ is monochromatic in $\Tilde{\alpha}$ (and in $\Tilde{\beta}$ by (i)), and
        \item[(v)] $s_\alpha' + s_\beta' \leq |X_\ell| + |Y_\ell|$.
    \end{itemize}
\end{claim}
\begin{proofclaim}
    Observe that in both $2$-dicolourings $\alpha'$ and $\beta'$, we are free to recolour any vertex of $X_\ell \cup Y_\ell$ since there is no monochromatic arc from $X$ to $Y$ and both $\vec{G}\ind{X}$ and $\vec{G}\ind{Y}$ are acyclic. Let $n_1$ (respectively $n_2$) be the number of vertices in $X_\ell \cup Y_\ell$ that are coloured 1 (respectively 2) in both $\alpha'$ and $\beta'$. Without loss of generality, assume that $n_1 \leq n_2$.
    Then we set each vertex of $X_\ell \cup Y_\ell$ to colour 2 in both $\alpha'$ and $\beta'$. Let $\Tilde{\alpha}$ and $\Tilde{\beta}$ the resulting $2$-dicolouring. Then $s_\alpha' + s_\beta'$ is exactly $|X_\ell| + |Y_\ell| + n_1 - n_2 \leq |X_\ell| + |Y_\ell|$.
\end{proofclaim}
\begin{claim}
    There is a redicolouring sequence between $\Tilde{\alpha}$ and $\Tilde{\beta}$ of length $|X_r^1| + |Y_r^1|$.
\end{claim}
\begin{proofclaim}
    By construction of $\Tilde{\alpha}$ and $\Tilde{\beta}$, we only have to exchange the colours of $x$ and $y$ for each arc $xy\in M_r^1$. Without loss of generality, we may assume that the colour of all vertices in $X_\ell \cup Y_\ell$ by $\Tilde{\alpha}$ and $\Tilde{\beta}$ is $1$. 
    
    We first prove that, by construction, we can recolour any vertex of $X_r^1\cup Y_r^1$ from 1 to 2. Assume not, then there is such a vertex $x\in X_r^1\cup Y_r^1$ such that recolouring $x$ from 1 to 2 creates a monochromatic directed cycle $C$. Since both $\vec{G}\ind{X}$ and $\vec{G}\ind{Y}$ are acyclic, $C$ must contain an arc of $M_r$. Since $M_r$ does not contain any monochromatic arc in $\Tilde{\alpha}$, then this arc must be incident to $x$. Now observe that colour 2, in $\Tilde{\alpha}$, induces an independent set on both $\vec{G}\ind{X}$ and $\vec{G}\ind{Y}$. This implies that $C$ must contain at least 2 arcs in $M_r$. This is a contradiction since recolouring $x$ creates exactly one monochromatic arc in $M_r$.
    
    Then, for each arc $xy\in M_r^1$, we can first recolour the vertex coloured 1 and then the vertex coloured 2. Note that we maintain the invariant that colour 2 induces an independent set on both $\vec{G}\ind{X}$ and $\vec{G}\ind{Y}$.
    We get a redicolouring sequence from $\Tilde{\alpha}$ to $\Tilde{\beta}$ in exactly $2|M_r^1| = |X_r^1| + |Y_r^1|$ steps.
\end{proofclaim}
Combining the three claims, we finally proved that there exists a redicolouring sequence between $\alpha$ and $\beta$ of length at most $n$.
\end{proof}

In the following, when $\alpha$ is a dicolouring of a digraph $D$, and $H$ is a subdigraph of $D$, we denote by $\alpha_{|H}$ the restriction of $\alpha$ to $H$.
We will prove Corollary~\ref{cor:Deltamin}, let us restate it.
\begingroup
    \def\thetheorem{\ref{cor:Deltamin}}
        \begin{corollary}
            Let $\vec{G}$ be an oriented graph of order $n$ with $\Delta_{\min}(\vec{G}) = \Delta \geq 1$, and let $k \geq \Delta + 1$. Then ${\cal D}_k(\vec{G})$ is connected and has diameter at most $2\Delta n$.
        \end{corollary}
    \addtocounter{theorem}{-1}
\endgroup
\begin{proof}
    
    We will show the result by induction on $\Delta$.
    
    Assume first that $\Delta=1$, let $k \geq 2$. Let $\alpha$ be any $k$-dicolouring of $\vec{G}$ and $\gamma$ be any 2-dicolouring of $\vec{G}$. To ensure that ${\cal D}_k(\vec{G})$ is connected and has diameter at most $2n$, it is sufficient to prove that there is a redicolouring sequence between $\alpha$ and $\gamma$ of length at most $n$.  Let $H$ be the digraph induced by the set of vertices coloured 1 or 2 in $\alpha$, and let $J$ be $V(\Vec{G}) \setminus V(H)$. By Theorem~\ref{thm:mindegree-1}, since $\Delta_{\min}(H) \leq \Delta_{\min}(\vec{G}) \leq 1$, we know that there exists a redicolouring sequence, in $H$, from $\alpha_{|H}$ to $\gamma_{|H}$ of length at most $|V(H)|$. This redicolouring sequence extends in $\vec{G}$ because it only uses colours 1 and 2. Let $\alpha'$ be the obtained dicolouring of $\vec{G}$. Since $\alpha'(v) = \gamma(v)$ for every $v\in H$, we can recolour each vertex in $J$ to its colour in $\gamma$. This shows that there is a redicolouring sequence between $\alpha$ and $\gamma$ of length at most $|V(H)| + |J| = |V(\vec{G})|$. This ends the case $\Delta= 1$.
    
    Assume now that $\Delta \geq 2$ and let $k \geq \Delta + 1$.
    Let $\alpha$ and $\beta$ be two $k$-dicolourings of $\vec{G}$. By Corollary~\ref{cor:brooks-oriented}, we know that $\dic(\vec{G}) \leq \Delta \leq k-1$.
    We first show that there is a redicolouring sequence of length at most $2n$ from $\alpha$ to some $(k-1)$-dicolouring  $\gamma$ of $\vec{G}$.
    From $\alpha$, whenever it is possible we recolour each vertex coloured $1,2$ or $k$ with a colour of $\{3,\dots,k-1\}$ (when $k=3$ we do nothing).  Let $\Tilde{\alpha}$ be the obtained dicolouring, and let $M$ be the set of vertices coloured in $\{3,\dots,k-1\}$ by $\Tilde{\alpha}$ (when $k=3$, $M$ is empty). We get that $H = \vec{G}-M$ satisfies $\Delta_{\min}(H) \leq 2$, since every vertex in $H$ has at least one in-neighbour and one out-neighbour coloured $c$ for every $c\in \{3,\dots,k-1\}$. By Corollary~\ref{cor:brooks-oriented}, there exists a 2-dicolouring $\gamma_{|H}$ of $H$. From $\Tilde{\alpha}_{|H}$, whenever it is possible, we recolour a vertex coloured 1 or 2 to colour $k$. Let $\hat{\alpha}$ be the resulting dicolouring, and $\hat{H}$ be the subdigraph of $H$ induced by the vertices coloured 1 or 2 in $\hat{\alpha}$. 
    We get that $\Delta_{\min}(\hat{H}) \leq 1$ since every vertex in $\hat{H}$ has, in $\Vec{G}$, at least one in-neighbour and one out-neighbour coloured $c$ for every $c\in \{3,\dots,k\}$.
    In at most $|V(\hat{H})|$ steps, using Theorem~\ref{thm:mindegree-1}, we can recolour the vertices of $V(\hat{H})$ to their colour in $\gamma_{|H}$ (using only colours 1 and 2).  
    Then we can recolour each vertex coloured $k$ to its colour in $\gamma_{|H}$.
    This results in a redicolouring sequence of length at most $2n$ from $\alpha$ to some $(k-1)$-dicolouring $\gamma$ of $\vec{G}$ , since colour $k$ is not used in the resulting dicolouring (recall that $M$ is coloured with $\{3,\dots,k-1\}$).
    
    Now, from $\beta$, whenever it is possible we recolour each vertex to colour $k$. Let $\Tilde{\beta}$ be the obtained $k$-dicolouring, and let $N$ be the set of vertices coloured $k$ in $\Tilde{\beta}$. We get that $J = \vec{G}-N$ satisfies $\Delta_{\min}(J) \leq \Delta - 1$. Thus, by induction, there exists a redicolouring sequence from $\Tilde{\beta}_{|J}$ to  $\gamma_{|J}$, in at most $2(\Delta-1)|V(J)|$ steps (using only colours $\{1,\dots,k-1\}$). Since $N$ is coloured $k$ in $\Tilde{\beta}$, this extends to a redicolouring sequence in $\vec{G}$. Now, since $\gamma$ does not use colour $k$, we can recolour each vertex in $N$ to its colour in $\gamma$.
    We finally get a redicolouring sequence from $\beta$ to $\gamma$ of length at most $2(\Delta-1)n$.
    Concatenating the redicolouring sequence from $\alpha$ to $\gamma$ and the one from $\gamma$ to $\beta$, we get a redicolouring sequence from $\alpha$ to $\beta$ in at most $2\Delta n$ steps.
\end{proof}
    \addtocounter{theorem}{4}
\section{An analogue of Brook's theorem for digraph redicolouring}
\label{section:recolouring_directed}

Let us restate Theorem~\ref{thm:extension_feghali}.

\begingroup
    \def\thetheorem{\ref{thm:extension_feghali}}
        \begin{theorem}
            Let $D$ be a connected digraph with $\Delta_{\max}(D) = \Delta \geq 3$, $k\geq \Delta+1$, and $\alpha$, $\beta$ two $k$-dicolourings of $D$. Then at least one of the following holds:
            \begin{itemize}
                \item $\alpha$ is $k$-frozen, or
                \item $\beta$ is $k$-frozen, or
                \item there is a redicolouring sequence of length at most $c_{\Delta}|V|^2$ between $\alpha$ and $\beta$, 
            where $c_{\Delta} = O(\Delta^2)$ is a constant depending only on $\Delta$.
            \end{itemize}
        \end{theorem}
    \addtocounter{theorem}{-1}
\endgroup

An \textit{$L$-assignment} of a digraph $D$ is a function which associates to every vertex a list of colours. An \textit{$L$-dicolouring} of $D$ is a dicolouring $\alpha$ where, for every vertex $v$ of $D$, $\alpha(v) \in L(v)$. An \textit{$L$-redicolouring sequence} is a redicolouring sequence $\gamma_1,\dots,\gamma_r$, such that for every $i\in \{1,\dots,r\}$, $\gamma_i$ is an $L$-dicolouring of $D$.

\begin{lemma}\label{lemma:blocked_vertex}
    Let $D=(V,A)$ be a digraph and $L$ be a list-assignment of $D$ such that, for every vertex $v\in V$, $|L(v)| \geq d_{\max}(v) +1$. Let $\alpha$ be an $L$-dicolouring of $D$. If $u\in V$ is blocked in $\alpha$, then for each colour $c\in L(u)$ different from $\alpha(u)$, $u$ has exactly one out-neighbour $u_c^+$ and one in-neighbour $u_c^-$  coloured $c$. Moreover, if $u_c^+ \neq u_c^-$, there must be a monochromatic directed path from $u_c^+$ to $u_c^-$. In particular, $u$ is not incident to a monochromatic arc.
\end{lemma}
\begin{proof}
    Since $u$ is blocked to its colour in $\alpha$, for each colour $c\in L(u)$ different from $\alpha(u)$, recolouring $u$ to $c$ must create a monochromatic directed cycle $C$. Let $v$ be the out-neighbour of $u$ in $C$ and $w$ be the in-neighbour of $u$ in $C$. Then $\alpha(v) = \alpha(w) = c$, and there is a monochromatic directed path (in $C$) from $v$ to $w$.
    
    This implies that, for each colour $c\in L(u)$ different from $\alpha(u)$, $u$ has at least one out-neighbour and at least one in-neighbour coloured $c$. Since $|L(u)| \geq d_{\max}(u) +1$, then $|L(u)| = d_{\max}(u) +1$, and $u$ must have exactly one out-neighbour and exactly one in-neighbour coloured $c$. In particular, $u$ cannot be incident to a monochromatic arc.
\end{proof}

\begin{lemma}\label{lemma:listdicolouring}
    Let $D=(V,A)$ be a digraph such that for every vertex $v\in V$, $N^+(v) \setminus N^-(v) \neq \emptyset$ and $N^-(v) \setminus N^+(v) \neq \emptyset$. Let $L$ be a list assignment of $D$, such that for every vertex $v\in V$, $|L(v)| \geq d_{\max}(v) +1$.
    
    Then for any pair of $L$-dicolourings $\alpha$, $\beta$ of $D$, there is an $L$-redicolouring sequence of length at most $(|V|+3) |V|$.
\end{lemma}

\begin{proof}
    Let $x = \diff (\alpha,\beta) = |\{v\in V \mid \alpha(v) \neq \beta(v) \}|$.
    We will show by induction on $x$ that there is an $L$-redicolouring sequence from $\alpha$ to $\beta$ of length at most $(|V|+3)x$.
    The result clearly holds for $x=0$ (i.e. $\alpha = \beta$).
    Let $v\in V$ be such that $\alpha(v) \neq \beta(v)$. We denote $\alpha(v)$ by $c$ and $\beta(v)$ by $c'$. If $v$ can be recoloured to $c'$, then we recolour it and we get the result by induction.
    
    Assume now that $v$ cannot be recoloured to $c'$. Whenever $v$ is contained in a directed cycle $C$ of length at least 3, such that every vertex of $C$ but $v$ is coloured $c'$, we do the following: we choose $w$ a vertex of $C$ different from $v$, such that $\beta(w) \neq c'$. We know that such a $w$ exists, for otherwise $C$ would be a monochromatic directed cycle in $\beta$. Now, since $w$ is incident to a monochromatic arc in $C$, and because $|L(w)| \geq d_{\max}(w)+1$, by Lemma~\ref{lemma:blocked_vertex}, we know that $w$ can be recoloured to some colour different from $c'$. Thus we recolour $w$ to this colour. Observe that it does not increase $x$.
    
    After repeating this process, maybe $v$ cannot be recoloured to $c'$ because it is adjacent by a digon to some vertices coloured $c'$. We know that these vertices are not coloured $c'$ in $\beta$. Thus, whenever such a vertex can be recoloured, we recolour it. After this, let $\eta$ be the obtained dicolouring. If $v$ can be recoloured to $c'$ in $\eta$, we are done. Otherwise, there must be some vertices, blocked to colour $c'$ in $\eta$, adjacent to $v$ by a digon. Let $S$ be the set of such vertices. Observe that, by Lemma~\ref{lemma:blocked_vertex}, for every vertex $s\in S$, $c$ belongs to $L(s)$, for otherwise $s$ would not be blocked in $\eta$. We distinguish two cases, depending on the size of $S$.
    \begin{itemize}
        \item If $|S|\geq 2$, then by Lemma~\ref{lemma:blocked_vertex}, $v$ can be recoloured to a colour $c''$, different from both $c$ and $c'$, because $v$ is adjacent by a digon with two neighbours coloured $c'$. 
        Hence we can successively recolour $v$ to $c''$,  and every vertex of $S$ to $c$ . This does not create any monochromatic directed cycle because for each $s\in S$, since $s$ is blocked in $\eta$, by Lemma~\ref{lemma:blocked_vertex} $v$ must be the only neighbour of $s$ coloured $c$ in $\eta$. 
        
        We can finally recolour $v$ to $c'$.
        \item If $|S|=1$, let $w$ be the only vertex in $S$. If $v$ can be recoloured to any colour (different from $c'$ since $w$ is coloured $c'$), then we first recolour $v$, allowing us to recolour $w$ to $c$, because $v$ is the single neighbour of $w$ coloured $c$ in $\eta$ by Lemma~\ref{lemma:blocked_vertex}. We finally can recolour $v$ to $c'$. 
        
        Assume then that $v$ is blocked to colour $c$ in $\eta$. Let us fix $w^+ \in N^+(w) \setminus N^-(w)$. Since $w$ is blocked to $c'$ in $\eta$, by Lemma~\ref{lemma:blocked_vertex}, there exists exactly one vertex $w^-\in N^-(w) \setminus N^+(w)$ such that $\eta(w^+) = \eta(w^-) = c''$ and there must be a monochromatic directed path from $w^+$ to $w^-$.
        
        Since $v$ is blocked to colour $c$ in $\eta$, either $vw^- \notin A$ or $w^+v \notin A$, otherwise, by Lemma~\ref{lemma:blocked_vertex}, there must be a monochromatic directed path from $w^-$ to $w^+$, which is blocking $v$ to its colour. But since there is also a monochromatic directed path from $w^+$ to $w^-$ (blocking $w$) there would be a monochromatic directed cycle, a contradiction (see Figure~\ref{fig:neighbourhood_of_w}).
        \begin{figure}[H]
            \begin{minipage}{\linewidth}
                \begin{center}	
                  \begin{tikzpicture}[thick,scale=1, every node/.style={transform shape}]
                    \tikzset{vertex/.style = {circle,fill=black,minimum size=6pt, inner sep=0pt}}
                    \tikzset{edge/.style = {->,> = latex'}}
            	 
                    \node[vertex, orange,label=above:$w$] (w) at  (0,0) {};
                    \node[vertex, g-green,label=below:$v$] (v) at  (0,-1) {};
                    \node[vertex, purple,label=right:$w^+$] (wp) at  (1.5,1) {};
                    \node[vertex, purple,label=left:$w^-$] (wm) at  (-1.5,1) {};
        
                    \draw[edge] (w) to (wp);
                    \draw[edge] (wm) to (w);
                    \draw [->,decorate,decoration=snake,purple] (wp) -- (wm);
                    \draw[edge, bend left = 20] (v) to (w);
                    \draw[edge, bend left = 20] (w) to (v);
                    \draw[edge, dashed] (v) to (wm);
                    \draw[edge, dashed] (wp) to (v);
                  \end{tikzpicture}
              \caption{The vertices $v,w,w^+$ and $w^-$.}
              \label{fig:neighbourhood_of_w}
            \end{center}    
          \end{minipage}
        \end{figure}
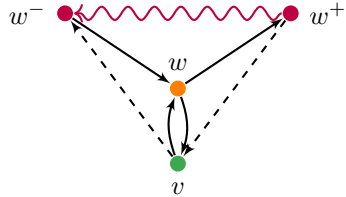
        We distinguish the two possible cases:
        \begin{itemize}
            \item if $vw^- \notin A$, then we start by recolouring $w^-$ with a colour that does not appear in its in-neighbourhood. This is possible because $w^-$ has a monochromatic entering arc, and because $|L(w^-)| \geq d_{\max}(w^-)+1$. We first recolour $w$ with $c''$, since $c''$ does not appear in its in-neighbourhood anymore ($w^-$ was the only one by Lemma~\ref{lemma:blocked_vertex}). Next we recolour $v$ with $c'$: this is possible because $v$ does not have any out-neighbour coloured $c'$ since $w$ was the only one by Lemma~\ref{lemma:blocked_vertex} and $w^-$ is not an out-neighbour of $v$. We can finally recolour $w$ to colour $c$ and $w^-$ to $c''$. After all these operations, we exchanged the colours of $v$ and $w$.
            
            \item if $w^+v \notin A$, then we use a symmetric argument.
        \end{itemize}
    \end{itemize}
    Observe that we found an $L$-redicolouring sequence from $\alpha$ to a $\alpha'$, in at most $|V|+3$ steps, such that $\diff(\alpha',\beta) < \diff(\alpha,\beta)$. Thus by induction, we get an $L$-redicolouring sequence of length at most $(|V|+3)x$ between $\alpha$ and $\beta$.
\end{proof}

We are now able to prove Theorem~\ref{thm:extension_feghali}. The idea of the proof is to divide the digraph $D$ into two parts. One of them is bidirected and we will use Theorem~\ref{thm:feghali} as a black box on it. In the other part, we know that each vertex is incident to at least two simple arcs, one leaving and one entering, and we will use Lemma~\ref{lemma:listdicolouring} on it.

\begin{proof} [Proof of Theorem~\ref{thm:extension_feghali}]
    Let $D=(V,A)$ be a connected digraph with $\Delta_{\max}(D) = \Delta$, $k\geq \Delta+1$. Let $\alpha$ and $\beta$ be two $k$-dicolourings of $D$.
    Assume that neither $\alpha$ nor $\beta$ is $k$-frozen.
    
    We first make a simple observation. For any simple arc $xy\in A$, we may assume that $N^+(y) \setminus N^-(y) \neq \emptyset$ and $N^-(x) \setminus N^+(x) \neq \emptyset$. If this is not the case, then every directed cycle containing $xy$ must contain a digon, implying that the $k$-dicolouring graph of $D$ is also the $k$-dicolouring graph of $D \setminus \{xy\}$. Then we may look for a redicolouring sequence in $D\setminus \{xy\}$.

    Let $X = \{v\in V \mid N^+(v) = N^-(v) \}$ and $Y = V \setminus X$. Observe that $D\ind{X}$ is bidirected, and thus the dicolourings of $D\ind{X}$ are exactly the colourings of $UG(D\ind{X})$. We first show that $\alpha_{|D\ind{X}}$ and $\beta_{|D\ind{X}}$ are not frozen $k$-colourings of $D\ind{X}$. If $Y$ is empty, then $D\ind{X} = D$ and $\alpha_{|D\ind{X}}$ and $\beta_{|D\ind{X}}$ are not $k$-frozen by assumption. Otherwise, since $D$ is connected, there exists $x\in X$ such that, in $D\ind{X}$, $d^+(x) = d^-(x) \leq \Delta - 1$, implying that $x$ is not blocked in any dicolouring of $D\ind{X}$.
    Thus, by Theorem~\ref{thm:feghali}, there is a redicolouring sequence $\gamma_1',\dots,\gamma_r'$ in $D\ind{X}$ from $\alpha_{|D\ind{X}}$ to $\beta_{|D\ind{X}}$, where $r \leq c_{\Delta}|X|^2$, and $c_{\Delta} = O(\Delta)$ is a constant depending on $\Delta$.
    
    We will show that, for each $i\in \{1,\dots, r-1\}$, if $\gamma_i$ is a $k$-dicolouring of $D$ which agrees with $\gamma_i'$ on $X$, then there exist a $k$-dicolouring $\gamma_{i+1}$ of $D$ that agrees with $\gamma_{i+1}'$ on $X$ and a redicolouring sequence from $\gamma_i$ to $\gamma_{i+1}$ of length at most $\Delta+2$.
    
    Observe that $\alpha$ agrees with $\gamma_1'$ on $X$. Now assume that there is such a $\gamma_i$, which agrees with $\gamma_i'$ on $X$, and let $v_i\in X$ be the vertex for which $\gamma_i'(v_i) \neq \gamma_{i+1}'(v_i)$. We denote by $c$ (respectively $c'$) the colour of $v_i$ in $\gamma_i'$ (respectively $\gamma_{i+1}'$). If recolouring $v_i$ to $c'$ in $\gamma_i$ is valid then we have the desired $\gamma_{i+1}$. Otherwise, we know that $v_i$ is adjacent with a digon (since $v_i$ is only adjacent to digons) to some vertices (at most $\Delta$) coloured $c'$ in $Y$. Whenever such a vertex can be recoloured to a colour different from $c'$, we recolour it. Let $\eta_i$ be the reached $k$-dicolouring after  these operations. If $v_i$ can be recoloured to $c'$ in $\eta_i$ we are done. If not, then the neighbours of $v_i$ coloured $c'$ in $Y$ are blocked to colour $c'$ in $\eta_i$. We denote by $S$ the set of these neighbours. We distinguish two cases:
    \begin{itemize}
        \item If $|S|\geq 2$, then by Lemma~\ref{lemma:blocked_vertex}, $v_i$ can be recoloured to a colour $c''$, different from both $c$ and $c'$, because $v_i$ has two neighbours with the same colour. Then we successively recolour $v_i$ to $c''$,  and every vertex of $S$ to $c$. This does not create any monochromatic directed cycle because, by Lemma~\ref{lemma:blocked_vertex}, for each $s\in S$, $v_i$ is the only neighbour of $s$ coloured $c$ in $\eta_i$. We can finally recolour $v_i$ to $c'$ to reach the desired $\gamma_{i+1}$.  
        \item If $|S|=1$, let $y$ be the only vertex in $S$. Since $y$ belongs to $Y$ and is blocked to its colour in $\eta_i$, by Lemma~\ref{lemma:blocked_vertex}, we know that $y$ has an out-neighbour $y^+ \in N^+(y) \setminus N^-(y)$ and an in-neighbour $y^- \in N^-(y) \setminus N^+(y)$ such that there is a monochromatic directed path from $y^+$ to $y^-$. Observe that both $y^+$ and $y^-$ are recolourable in $\eta_i$ by Lemma~\ref{lemma:blocked_vertex}, because there are incident to a monochromatic arc.
        
        \begin{itemize}
            \item If $v_i$ is not adjacent to $y^+$, then we recolour $y^+$ to any possible colour, and we recolour $y$ to $\eta_i(y^+)$. We can finally recolour $v_i$ to $c'$ to reach the desired $\gamma_{i+1}$. 
            
            \item If $v_i$ is not adjacent to $y^-$, then we recolour $y^-$ to any possible colour, and we recolour $y$ to $\eta_i(y^-)$. We can finally recolour $v_i$ to $c'$ to reach the desired $\gamma_{i+1}$.
            
            \item Finally if $v_i$ is adjacent to both $y^+$ and $y^-$, since $\eta_i(y^+) = \eta_i(y^-)$, then $v_i$ can be recoloured to a colour $c''$ different from $c$ and $c'$. This allows us to recolour $y$ to $c$, and we finally can recolour $v_i$ to $c'$ to reach the desired $\gamma_{i+1}$. 
        \end{itemize}
    \end{itemize}
    
    We have shown that there is a redicolouring sequence of length at most $(\Delta+2)c_{\Delta}n^2$ from $\alpha$ to some $\alpha'$ that agrees with $\beta$ on $X$.
    Now we define the list-assignment: for each $y\in Y$, 
    \[L(y) = \{1,\dots,k\} \setminus \{ \beta(x) \mid x \in N(y) \cap X \}. \] 
    Observe that, for every $y\in Y$, 
    \[|L(y)| \geq k - |N^+(y) \cap X| \geq \Delta + 1 - (\Delta - d^+_Y(y))  \geq d^+_Y(y) + 1.\]
    Symmetrically, we get $|L(y)| \geq d^-_Y(y) +1 $. This implies, in $D\ind{Y}$, $|L(y)| \geq d_{\max}(y) + 1$. Note also that both $\alpha'_{|D\ind{Y}}$ and $\beta_{|D\ind{Y}}$ are $L$-dicolourings of $D\ind{Y}$. Note finally that, for each $y\in Y$, $N^+(y) \setminus N^-(y) \neq \emptyset$ and $N^+(y) \setminus N^-(y) \neq \emptyset$ by choice of $X$ and $Y$ and by the initial observation.
    By Lemma~\ref{lemma:listdicolouring}, there is an $L$-redicolouring sequence in $D\ind{Y}$ between $\alpha'_{|D\ind{Y}}$ and $\beta_{|D\ind{Y}}$, with length at most  $(|Y|+3)|Y|$. By choice of $L$, this extends directly to a redicolouring sequence from $\alpha'$ to $\beta$ on $D$ of the same length.
    
    The concatenation of the redicolouring sequence from $\alpha$ to $\alpha'$ and the one from $\alpha'$ to $\beta$ leads to a redicolouring sequence from $\alpha$ to $\beta$ of length at most $c_{\Delta}'|V|^2$, where $c_\Delta' =O(\Delta^2)$ is a constant depending on $\Delta$.
\end{proof}

\begin{remark}
    If $\alpha$ is a $k$-frozen dicolouring of a digraph $D$, with $k\geq \Delta_{\max}(D) +1$, then $D$ must be bidirected. If $D$ is not bidirected, then we choose $v$ a vertex incident to a simple arc. If $v$ cannot be recoloured in $\alpha$, by Lemma~\ref{lemma:blocked_vertex}, since $v$ is incident to a simple arc, there exists a colour $c$ for which $v$ has an out-neighbour $w$ and an in-neighbour $u$ both coloured $c$, such that $u\neq w$ and there is a monochromatic directed path from $w$ to $u$. But then, every vertex on this path is incident to a monochromatic arc, and it can be recoloured by Lemma~\ref{lemma:blocked_vertex}. Thus, $\alpha$ is not $k$-frozen. This shows that an obstruction of Theorem~\ref{thm:extension_feghali} is exactly the bidirected graph of an obstruction of Theorem~\ref{thm:feghali}.
\end{remark}
\section{Further research}
\label{section:open_problems}

In this paper, we established some analogues of Brooks' Theorem for the dichromatic number of oriented graphs and for digraph redicolouring. Many open questions arise, we detail a few of them.

Restricted to oriented graphs, Mcdiarmid and Mohar (see~\cite{harutyunyanEJC18}) conjectured that the Directed Brooks' Theorem can be improved to the following.
\begin{conjecture}[Mcdiarmid and Mohar]
    Every oriented graph $\vec{G}$ has $\dic(\vec{G}) = O\left(\frac{\Delta_{\max}}{\log (\Delta_{\max})}\right)$.
\end{conjecture}

\medskip

Concerning digraph redicolouring, we believe that Corollary~\ref{cor:Deltamin} and Theorem~\ref{thm:extension_feghali} can be improved. We pose the following conjecture.
\begin{conjecture}
There is an absolute constant $c$ such that 
for  every integer $k$ and every oriented graph $\vec{G}$ on $n$ vertices, such that $k \geq \Delta_{\min}(\vec{G}) + 1$, the diameter of ${\cal D}_k(\vec{G})$ is bounded by $c n$.
\end{conjecture}

\begin{conjecture}
There is an absolute constant $d$ such that 
for  every integer $k$ and every digraph $D$ on $n$ vertices, with $k \geq \Delta_{\max}(D) + 1 \geq 4$, the diameter of ${\cal D}_k(D)$ is bounded by $d n^2$.
\end{conjecture}

It would also be nice to extend Theorem~\ref{thm:bousquet_al} to directed graphs, that is to show that the diameter of ${\cal D}_k(D)$ is bounded by $f(\Delta_{\max}(D)) n$ for some computable function $f$ whenever $k \geq \Delta_{\max}(D) + 1 \geq 4$. To prove it, it would be sufficient to show the analogue of Lemma~\ref{lemma:listdicolouring} with an $L$-redicolouring sequence of length at most $f(\Delta_{\max}(D))|V|$, and then follow the proof of Theorem~\ref{thm:extension_feghali} (using Theorem~\ref{thm:bousquet_al} instead of Theorem~\ref{thm:feghali}).

\medskip 

Given an orientation $\vec{G}$ of a planar graph, a celebrated conjecture from Neumann-Lara~\cite{neumannlaraJCT33} states that the dichromatic number of $\vec{G}$ is at most 2. It is known that it must be 4-mixing because planar graphs are 5-degenerate~\cite{papierAMADEUS}. It is also known that there exists 2-freezable orientations of planar graphs~\cite{papierAMADEUS}. Thus the following problem, stated in~\cite{papierAMADEUS}, remains open:
\begin{question}
    Is every oriented planar graph 3-mixing ?
\end{question}

\section*{Acknowledgement}
I am grateful to Frédéric Havet and Nicolas Nisse for stimulating discussions.

\bibliography{refs}

\end{document}